\documentclass[10pt,twoside]{article}
\usepackage[latin1]{inputenc}
\usepackage{amsmath}
\usepackage{dsfont}
\usepackage{graphicx}
\usepackage{yhmath}
\usepackage[table]{xcolor}
\usepackage{mathrsfs} 
\usepackage{amssymb}
\usepackage{url}
\usepackage{makecell}
\usepackage{arydshln}
\usepackage{multirow}
\usepackage{arydshln}
\usepackage{mathtools}
\usepackage{stmaryrd}  

\input epsf
\def\limiten{\renewcommand{\arraystretch}{0.5}
\begin{array}[t]{c}\stackrel{}{\longrightarrow} \\
{\scriptstyle n\rightarrow
\infty}\end{array}\renewcommand{\arraystretch}{1}}

\def\limitepsn{\renewcommand{\arraystretch}{0.5}
\begin{array}[t]{c}\stackrel{a.s.}{\longrightarrow} \\
{\scriptstyle n \rightarrow
\infty}\end{array}\renewcommand{\arraystretch}{1}}

\def\limiteloin{\renewcommand{\arraystretch}{0.5}
\begin{array}[t]{c}\stackrel{{\cal D}}{\longrightarrow} \\
{\scriptstyle n\rightarrow
\infty}\end{array}\renewcommand{\arraystretch}{1}}

\def\limiteproban{\renewcommand{\arraystretch}{0.5}
\begin{array}[t]{c}\stackrel{{\cal P}}{\longrightarrow} \\
{\scriptstyle n\rightarrow
\infty}\end{array}\renewcommand{\arraystretch}{1}}

\numberwithin{equation}{section}

\newtheorem{thm}{Theorem}[section]

\newtheorem{lem}[thm]{Lemma}

\newcommand{\E}{\ensuremath{\mathbb{E}}}
\newcommand{\R}{\ensuremath{\mathbb{R}}}
\newcommand{\Z}{\ensuremath{\mathbb{Z}}}

\newcommand{\N}{\ensuremath{\mathbb{N}}}

\definecolor{grisclair}{gray}{0.9}

\renewcommand{\arraystretch}{.8}

\begin{document}
\title{\bf Epidemic change-point detection in  general causal
time series}
 \maketitle \vspace{-1.0cm}
\begin{center}
   Mamadou Lamine DIOP \footnote{Supported by
   the MME-DII center of excellence (ANR-11-LABEX-0023-01) 
   } 
   and 
     William KENGNE \footnote{Developed within the ANR BREAKRISK: ANR-17-CE26-0001-01 and the  CY Initiative of Excellence (grant "Investissements d'Avenir" ANR-16-IDEX-0008), Project "EcoDep" PSI-AAP2020-0000000013} 
 \end{center}

  \begin{center}
  { \it 
 THEMA, CY Cergy Paris Université, 33 Boulevard du Port, 95011 Cergy-Pontoise Cedex, France\\
  E-mail: mamadou-lamine.diop@u-cergy.fr ; william.kengne@u-cergy.fr  \\
  }
\end{center}

 \pagestyle{myheadings}
 \markboth{Epidemic change-point detection in  general causal
time series}{Diop and Kengne}

~~\\
\textbf{Abstract}:
 We consider an epidemic change-point detection in a large class of causal time series models, including among other processes, AR($\infty$), ARCH($\infty$), TARCH($\infty$), ARMA-GARCH.
A test statistic based on the Gaussian quasi-maximum likelihood estimator of the parameter is proposed.
%
It is shown that, under the null hypothesis of no change, the test statistic converges to a distribution obtained from a difference of two Brownian bridge and diverges to infinity under the epidemic alternative.
Numerical results for simulation and real data example are provided.

 \medskip
 

\section{Introduction}
%
%
%

We consider a general class of affine causal  time series models  in a semiparametric setting.  
 Let $(\xi_t)_{t\in \Z}$ be a sequence of centered independent
 and identically distributed (iid)  random variables satisfying $\E \xi_0^2= 1$  and $\Theta $ a compact subset of $\R^{d}$ ($d \in \N$).
 For $ \mathcal T \subset \Z$ and any $ \theta \in \Theta$, define ~\\
  \textbf{Class} $\mathcal{AC}_{\mathcal T}(M_{\theta},f_{\theta})$: A process $\{X_{t},\,t\in \mathcal T \}$ belongs to $\mathcal{AC}_{\mathcal T}(M_\theta,f_\theta)$ if it satisfies:
   \begin{equation}\label{Model} 
     X_t =M_{\theta}(X_{t-1}, X_{t-2}, \ldots) \cdot \xi_t + 
      f_{\theta}(X_{t-1}, X_{t-2}, \ldots)~~\forall t \in \mathcal T,
   \end{equation}
 where $M_\theta , f_\theta : \R^{\infty} \rightarrow \R$ are two measurable functions.
 The existence of a stationary and ergodic solution as well as the inference   for the class $\mathcal{AC}_{\mathcal \Z}(M_{\theta^*},f_{\theta^*})$ have been addressed by \nocite{Bardet2009} Bardet and Wintenberger (2009). 
 Numerous classical time series such as AR($\infty$), ARCH($\infty$), TARCH($\infty$) or ARMA-GARCH models belong to this class (see Bardet and Wintenberger (2009)). 
 This class of models has now been well, see for instance \nocite{Bardet2012} Bardet {\it et al.} (2012), \nocite{Kengne2012} Kengne (2012), \nocite{Bardet2014} Bardet and Kengne (2014) for change-point detection on this class; \nocite{Bardet2017} Bardet {\it et al.} (2017) for inference based on  the Laplacian quasi-likelihood; \nocite{Bardet2020} Bardet {\it et al.} (2020), \nocite{Kengne2021strongly} Kengne (2020) for model selection in this class. 
 
 \medskip

 We focus here on the epidemic change-point detection in the class $\mathcal{AC}_{\mathcal T}(M_{\theta},f_{\theta})$. 
 Assume that a trajectory $(X_{1},\cdots,X_{n})$ of the process $\{X_{t},\,t\in \Z \}$ is observed and consider the following test hypotheses:
\begin{enumerate}
    \item [ H$_0$:] $(X_1,\cdots,X_n)$  is a trajectory of the process $\{X_{t},\,t\in \Z \} \in \mathcal{AC}_{\mathcal \Z}(M_{\theta_0^*},f_{\theta_0^*})$  with $\theta_0^* \in \Theta$.
 
\item[ H$_1$:]: there exists $(\theta^{*}_1,\theta^{*}_2,t^{*}_1,t^{*}_2) \in \Theta^{2}\times \{2,3,\cdots, n-1 \}^2$ (with $\theta^{*}_1 \neq \theta^{*}_2$ and $t^{*}_1 <t^{*}_2$) such that $(X_1,\cdots,X_n)$ belongs to $ \mathcal{AC}_{ \{1,\cdots,t_1^* \} }(M_{\theta_1^*},f_{\theta_1^*}) \bigcap \mathcal{AC}_{ \{t^*_1+1,\cdots t^*_2 \} }(M_{\theta_2^*},f_{\theta_2^*}) \bigcap \mathcal{AC}_{ \{t^*_2+1,\cdots n \} }(M_{\theta_1^*},f_{\theta_1^*}) $.
\end{enumerate} 
 The epidemic alternative H$_1$ refers to the so-called epidemic period, which runs from $t^*_1$ to $t^*_2$. 
 
  \medskip
  
 Several works in the literature are devoted to the epidemic change-point detection in time series.
We refer among others, to \nocite{Levin1985} Levin and Kline (1985), \nocite{Yao1993} Yao (1993), \nocite{Csorgo1997} Cs{\"o}rg{\"o} and Horv{\'a}th (1997), \nocite{Ramanayake2003} Ramanayake and Gupta (2003), \nocite{Ravckauskas2004} Ra{\v{c}}kauskas and Suquet (2004), \nocite{Ravckauskas2006} Ra{\v{c}}kauskas and Suquet (2006), \nocite{Guan2007} Guan (2007), \nocite{Jaruvskova2011log} Jaru{\v{s}}kov{\'a} and Piterbarg (2011), \nocite{Aston2012detecting, Aston2012Evaluating} Aston and Kirch (2012a, 2012b),  \nocite{Bucchia2014} Bucchia (2014), \nocite{Graiche2016} Graiche {\it et al.} (2016).
As pointed out by \nocite{Diop2021epidemic} Diop and Kengne (2021), most of these procedures are developed for the epidemic change-point detection in the mean of random variables. 
The latter authors addressed this issue for a general class of integer valued time series.

\medskip

In this new contribution, we propose a test based on the Gaussian quasi-likelihood for the epidemic change-point detection in the class of  affine causal models $\mathcal{AC}_{\mathcal T}(M_{\theta},f_{\theta})$. 
Under the null hypothesis of no change, the proposed statistic converges to a distribution obtained from a difference between two Brownian bridges; this statistic diverges to infinity under the epidemic alternative.
These findings lead to a test which has correct size asymptotically and is consistent in power.  

 \medskip

The rest of the paper is outlined as follows. Section 2 provides some assumptions and the definition of the Gaussian quasi-likelihood.
 Section 3 focuses on the construction of the test statistic and
the asymptotic studies under the null and the epidemic alternative. 
Some numerical results for simulation and real data example are displayed in
Section 4. Section 5 is devoted to the proofs of the main results.

\section{Assumptions and  QMLE}\label{Sect_Ass_PQMLE}
Throughout the sequel, we use the following notations:
{\em
\begin{itemize}
 \item $ \|x \| \coloneqq  \sqrt{\sum_{i=1}^{p} |x_i|^2 } $, for any $x \in \mathbb{R}^{p}$;
%
%
 \item $ \|x \| \coloneqq  \sqrt{\sum_{i=1}^{p} \sum_{j=1}^{q} |x_{i,j}|^2 } $, for any matrix $x=(x_{i,j}) \in M_{p,q}(\R)$; where $M_{p,q}(\R)$ denotes the set of matrices of dimension $p\times q$ with coefficients in $\R$;
\item  $\left\|f\right\|_{\Theta} \coloneqq \sup_{\theta \in \Theta}\left(\left\|f(\theta)\right\|\right)$ for any function $f:\Theta \longrightarrow   M_{p,q}(\R)$;
\item $\left\|Y\right\|_r \coloneqq \E\left(\left\|Y\right\|^r\right)^{1/r}$, where $Y$ is a random vector with finite $r-$order moments; 
\item $T_{\ell,\ell'}=\{\ell,\ell+1,\cdots,\ell'\}$ for any $(\ell,\ell') \in \N^2$ such as $\ell \leq \ell'$.
\end{itemize}
}

 \medskip
    \noindent 
In the sequel, 0 denote the null vector of any vector space. 
%
%
\noindent For $\Psi_\theta=f_\theta, \, M_\theta$ and any compact set $ \mathcal{K}\subseteq \Theta$, define

\medskip

{\bf Assumption A$_i(\Psi_\theta, \mathcal{K})$} ($i=0,1,2$): {\em Assume that
$\|{\partial^i\Psi_\theta(0)}/{\partial\theta^i}\|_{\mathcal{K}}<\infty$
 and there exists a sequence of non-negative real number $(\alpha^{(k)}_i(\Psi_\theta, \mathcal{K}))_{i\geq 1}$ such that $ \sum\limits_{k=1}^{\infty} \alpha^{(i)}_k(\Psi_\theta, \mathcal{K}) <\infty$  satisfying
\begin{equation*}
\Big\|\dfrac{\partial^i\Psi_\theta(x)}{\partial\theta^i}-\dfrac{\partial^i\Psi_\theta(y)}{\partial\theta^i}\Big\|_{\mathcal{K}}
\leq \sum\limits_{k=1}^{\infty}\alpha^{(i)}_k(\Psi_\theta, \mathcal{K} )|x_k-y_k|\quad \text{for all}~x, y \in \R^{\infty},
\end{equation*}}
where $x$, $y$, $x_k$, $y_k$ are respectively replaced by $x^2$, $y^2$, $x_k^2$, $y_k^2$ if $\Psi_\theta = h_\theta := M_\theta^2$.

\medskip

   For any $r\geq 1$, define 
\begin{multline*}\label{Set_Theta(r)}
\Theta(r) = \big\{
 \theta \in \R^d \, \big / \, \textbf{A}_0 (f_\theta,\{\theta\}) \   \text{and}\ \textbf{A}_0 (M_\theta,\{\theta\})   \    \text{hold with} 
 ~ \sum\limits_{k=1}^{\infty} \left\{\alpha^{(0)}_{k}(f_\theta,\{\theta\}) + \|\xi_0\|_r \alpha^{(0)}_{k}(M_\theta,\{\theta\}\right\} 
 <1
\big\}\\
\bigcup 
\big\{
 \theta  \in \R^d \ \big / \ f_\theta=0 \text{ and } \textbf{A}_0 (h_\theta,\{\theta\})  \text{ holds with } 
 \|\xi_0\|^2_r  \sum\limits_{k=1}^{\infty} \alpha^{(0)}_{k}(h_\theta,\{\theta\}) 
 <1
\big\}.
\end{multline*}
These Lipschitz-type conditions are notably useful when studying the existence of solutions of the class $\mathcal{AC}_{\mathcal T}(M_{\theta},f_{\theta})$. 
 If $\theta\in\Theta(r)$, then there exists a $\tau$-weakly dependent stationary and ergodic solution   $X=(X_t)_{t\in \Z}\in  \mathcal{AC}_{\Z}(M_{\theta},f_{\theta})$ satisfying $\| X_0 \|_r <\infty$ (see \nocite{Doukhan2008} Doukhan and Wintenberger (2008) and  \nocite{Bardet2009} Bardet and Wintenberger (2009)).

\medskip

Consider a trajectory $(X_1,\cdots,X_n)$ of a process $X=(X_t)_{t \in \Z}$. If $(X_1,\cdots,X_n) \in \mathcal{AC}_{ \{1,\cdots,n \}  }(M_{\theta},f_{\theta}) $, then for any segment $ T \subset \{1,\cdots,n \} $, the conditional Gaussian quasi-(log)likelihood computed on $T$ is given by,
\begin{equation}\label{defML}
L(T,\theta) := -\dfrac{1}{2} \sum\limits_{t \in T} q_t(\theta)   ~~ \text{with}  ~~ q_t(\theta)= \dfrac{(X_t-f^t_{\theta})^2}{h_\theta^t} + \log(h_\theta^t) 
\end{equation} 
       where $ f^t_{\theta}=f_{\theta}\big(X_{t-1},X_{t-2}\ldots\big) $, $ M^t_{\theta}=M_{\theta}\big(X_{t-1},X_{t-2}\ldots\big) $ and
   $ h^t_{\theta}= ({M^{t}_{\theta}})^2 $. 
  In the sequel, we deal with an approximated quasi-(log)likelihood contrast given for any segment  $ T \subset \{1,\cdots,n \} $ by,
 \[ \widehat L(T,\theta):=-\frac{1}{2}\sum\limits_{t \in T}\widehat q_t(\theta)\quad
 \textrm{where}\;\;\;\widehat{q}_t(\theta):=\frac{\big(X_t-\widehat{f}^t_{\theta}\big)^2}{\widehat{h}{^t_{\theta}}} +\log\big(\widehat{h}{^t_{\theta}}\big) \]
 with
 $ \widehat{f}^t_{\theta}=f_{\theta}\big(X_{t-1},\ldots,X_{1},0,0,\cdots\big)$, $ \widehat{M}^t_{\theta}=
 M_{\theta}\big(X_{t-1},\ldots,X_{1},0,0,\cdots\big) $ and $\widehat{h}^t_{\theta}=( \widehat{M}^{t}_{\theta} )^2$; and consider the estimator,
\begin{equation}\label{QMLE}
\widehat{\theta}_n(T) :=  \underset{\theta\in \Theta} {\text{argmax}}(\widehat{L}(T,\theta))
\end{equation}
The following assumptions are needed to study the asymptotic behavior of the estimator defined in (\ref{QMLE}).
 \noindent {\bf Assumption D$(\Theta)$:} $\exists\underline{h}>0$ such that
$\displaystyle \inf_{\theta \in
 \Theta}(h_\theta(x))\geq \underline{h}$ for all $x\in \R^{\N}.$ 
 
\medskip
 \noindent {\bf Assumption Id($\Theta$):} For a process $(X_t)_{t \in \Z} \in \mathcal{AC}_{\Z}(M_{\theta^*},f_{\theta^*})$ and for all  $\theta \in \Theta$,
\[ \Big( f_{\theta^*}(X_0,X_{-1},\cdots)=f_{\theta}(X_0,X_{-1},\cdots)~\text{and}~h_{\theta^*}(X_0,X_{-1},\cdots)=h_{\theta}(X_0,X_{-1},\cdots) \ \text{a.s.}\Big) \Rightarrow \ \theta^* = \theta. \]
{\bf Assumption Var($\Theta$):} For a process $(X_t)_{t \in \Z} \in \mathcal{AC}_{\Z}(M_{\theta^*},f_{\theta^*})$, one
 of the families $ \big( \dfrac{\partial f_{\theta^*}}{\partial \theta^{i}}(X_0,X_{-1},\cdots) \big)_{1\leq i \leq d}$ or 
 $\big( \dfrac{\partial h_{\theta^*}}{\partial \theta^{i}}(X_0,X_{-1},\cdots) \big)_{1\leq i \leq d}  \quad $
 is $a.e.$ linearly independent.

\medskip

\medskip

\noindent Under H$_0$ and the above assumptions, Bardet and Wintenberger (2009) established the consistency and the asymptotic normality of the estimator $\widehat{\theta}_n(T_{1,n}) $ for the class $\mathcal{AC}_{\Z}(M_{\theta^*_0},f_{\theta^*_0})$.

\section{Test statistic and asymptotic results}
 Under H$_0$, recall that (see \nocite{Bardet2009}  Bardet and Wintenberger (2009)), for the class $\mathcal{AC}_{\Z}(M_{\theta^*_0},f_{\theta^*_0})$, it holds that
 \begin{equation}\label{TLC}
 \sqrt n \big (\widehat{\theta}(T_{1,n})-\theta_0^* \big) \limiteloin {\cal N}\big ( 0\, , \,{F}^{-1}\,G\,{F}^{-1} \big ),
 \end{equation}
with
\begin{equation}\label{FG}
 \displaystyle G:= \E \Big [ \dfrac{\partial q_0(\theta^*_0)}{\partial \theta} \dfrac{\partial q_0(\theta^*_0)}{\partial \theta} ' \Big] \quad \mbox {and}\quad \displaystyle  F:=  \E \Big [ \dfrac{\partial^2 q_0(\theta^*_0)}{\partial \theta \partial \theta'}  \Big],
  \end{equation}
 where $'$ denotes the transpose.
 For any segment $ T \subset \{1,\cdots,n \} $, consider the following matrices,
 \begin{equation}\label{hatFG}
 \widehat{G}(T):= \dfrac{1}{\mbox{Card}(T)} \sum \limits_{t \in T} \Big( \dfrac{\partial \widehat{q}_{t}( \widehat{\theta}(T))}{\partial \theta} \Big)
                       \Big( \dfrac{\partial \widehat{q}_{t}( \widehat{\theta}(T))}{\partial \theta} \Big)'
                       \quad\mbox{and}\quad
    \widehat{F}(T):=  \dfrac{1}{\mbox{Card}(T)} \sum \limits_{t \in T} \dfrac{\partial^2 \widehat{q}_t(\widehat{\theta}(T))}{\partial \theta \partial \theta'}.
   \end{equation}
    Under H$_0$, $\widehat{G}(T_{1,n})$ and $\widehat{F}(T_{1,n})$  are consistent estimators of $G$ and $F$, respectively.   

\medskip
\noindent
In the sequel, we follow the idea of \nocite{Diop2021EpidCount} Diop and Kengne (2021). 
Let $(u_n)_{n\geq 1}$, $(v_n)_{n\geq 1}$ be two integer valued sequences such that:  $ (u_n,v_n) =o(n)$ and $ u_n,v_n \limiten +\infty$.
For all $n \geq 1$, define the matrix 
\begin{multline*}\label{Sigma_un}
\widehat{\Sigma}(u_n)=\frac{1}{3} 
\big[
\widehat F(T_{1,u_n})  \widehat G(T_{1,u_n})^{-1}   \widehat F(T_{1,u_n}) 
+
\widehat F(T_{u_n+1,n-u_n})  \widehat G(T_{u_n+1,n-u_n})^{-1}  \widehat F(T_{u_n+1,n-u_n}) \\
+
\widehat F(T_{n-u_n+1,n})  \widehat G(T_{n-u_n+1,n})^{-1}  \widehat F(T_{n-u_n+1,n})
\big]
\end{multline*}
where $\widehat G(T_{1,u_n})^{-1}$, $\widehat G(T_{u_n+1,n-u_n})^{-1}$, $\widehat G(T_{n-u_n+1,n})^{-1}$ are replaced by 0 if these matrices are not invertible.
Also, define the set
\[
\mathcal T_n = \left\{ (k_1, k_2) \in {([v_n, n-v_n] \cap \N)}^2 ~\text{ with }~  k_2- k_1  \geq v_n \right\}.
\]
For all $(k_1,k_2) \in \mathcal  T_n$, set
 \begin{equation}\label{def_stat_C_n,k1,k2}
C_{n,k_1,k_2}= \frac{(k_2-k_1)}{n^{3/2}} \left[ \left(n-(k_2-k_1)\right)\widehat{\theta}(T_{k_1+1,k_2}) -  k_1 \widehat{\theta}(T_{1,k_1}) - (n-k_2)\widehat{\theta}(T_{k_2+1,n}) \right],
\end{equation}
and consider the test statistic
 \begin{equation}\label{Stat_Q_n}
 \widehat{Q}_n=\max_{(k_1, k_2) \in \mathcal T_n} \widehat{Q}_{n,k_1,k_2}
~\text{ with }~
 \widehat{Q}_{n,k_1,k_2}=C_{n,k_1,k_2}' \widehat{\Sigma}(u_n) C_{n,k_1,k_2}.
  \end{equation}
As pointed out by \nocite{Diop2021EpidCount} Diop and Kengne (2021), this test statistic coincides with those proposed by 
\nocite{Ravckauskas2004holder} Rackauskas and Suquet (2004) (statistic $UI(n,\rho)$), \nocite{Jaruvskova2011log} Jarusková and Piterbarg (2011) (statistic $T_1^2$), \nocite{Bucchia2014testing} Bucchia (2014) (statistic $T_n(\alpha,\beta)$) or \nocite{Aston2012detecting} Aston and Kirch (2012) (statistic $T_n^{B_2}$)
 for the particular case of epidemic change-point detection in the mean. In this sense, the test considered here can be seen as a generalization these procedures.  

\medskip
\noindent
 The following theorem provides the asymptotic behavior of the statistic $\widehat{Q}_n$ under the null hypothesis. In the condition  (\ref{eq_th1}) in this theorem, we make the convention that if \textbf{A$_i(M_\theta,\Theta)$} holds, then $\alpha^{(i)}_{k}(h_\theta,\Theta) = 0$ for all $k \in \N$ and if \textbf{A$_i(h_\theta,\Theta)$} holds, then $\alpha^{(i)}_{k}(M_\theta,\Theta) = 0$ for all $k \in \N$.    
\begin{thm}\label{th1}
Under H$_0$ with $\theta^*_0 \in \overset{\circ}{\Theta} \cap \Theta(4)$, assume that \textbf{D$(\Theta)$}, \textbf{Id($\Theta$)}, \textbf{Var($\Theta$)} (for the class $\mathcal{AC}_{\Z}(M_{\theta^*_0},f_{\theta^*_0})$),  \textbf{A$_i(f_\theta,\Theta)$}, \textbf{A$_i(M_\theta,\Theta)$} (or \textbf{A$_i(h_\theta,\Theta)$}) hold with
\begin{equation}\label{eq_th1}
 \alpha^{(i)}_{k}(f_\theta,\Theta)+\alpha^{(i)}_{k}(M_\theta,\Theta)
 +\alpha^{(i)}_{k}(h_\theta,\Theta) = O(k^{-\gamma}) 
 \text{ for } i=0,1,2 \text{ and some } \gamma >3/2. 
\end{equation}

Then,
\begin{equation}\label{res_th1}
\widehat{Q}_n \limiteloin \sup_{0\leq \tau_1 < \tau_2 \leq 1}\left\|W_d(\tau_1)-W_d(\tau_2)\right\|^{2}, 
\end{equation}
where $W_d$ is a $d$-dimensional Brownian bridge.
\end{thm}

\noindent 
For any $\alpha \in (0,1)$, denote $c_{d,\alpha}$ the $(1-\alpha)$-quantile of the distribution
of $\sup_{0\leq\tau_1 < \tau_2\leq 1}\left\|W_d(\tau_1)-W_d(\tau_2)\right\|^{2}$.
Therefore, at a nominal level $\alpha \in (0,1)$, the critical region of the test is $(\widehat{Q}_{n}>c_{d,\alpha})$; 
 which leads to a procedure with correct size asymptotically. 
Table 1 of Diop and Kengne (2021) provides the values of $c_{d,\alpha}$ for $\alpha=0.01,\,0.05,\,0.10$ and $d=1,\ldots,5$.

\medskip

\noindent
For asymptotic under the epidemic alternative, the following additional condition is needed.

  \medskip
 \noindent  {\bf Assumption B}: {\em There exists $(\tau^*_1,\tau^*_2) \in (0,1)^2$ such that $(t^*_1,t^*_2)=([n\tau^*_1], [n\tau^*_2])$ (with $[\cdot]$ is the integer part).}

\medskip

\noindent
We have the following result.
\begin{thm}\label{th2}
Under $H_1$ with $\theta^*_1, \theta^*_2 \in  \overset{\circ}{\Theta} \cap \Theta(4)$, assume that \textbf{D$(\Theta)$}, \textbf{Id($\Theta$)}, \textbf{Var($\Theta$)} (for the classes $\mathcal{AC}_{\Z}(M_{\theta^*_1},f_{\theta^*_1})$ and $\mathcal{AC}_{\Z}(M_{\theta^*_2},f_{\theta^*_2})$), \textbf{A$_i(f_\theta,\Theta)$}, \textbf{A$_i(M_\theta,\Theta)$} (or \textbf{A$_i(h_\theta,\Theta)$}) and (\ref{eq_th1}) hold. 
Then,

\begin{equation}\label{res_th2}
 \widehat{Q}_n  \limiteproban +\infty .
 \end{equation}
\end{thm}
 This theorem shows that the proposed procedure is consistency in power. 
 An estimator of the change-points $\underline t^*=(t^*_1,t^*_2)$ under the epidemic alternative is given by 
 \[   \widehat{ \underline t}_n =  \underset{(k_1, k_2) \in \mathcal T_n}{ \text{argmax}  }  C_{n,k_1,k_2}' \widehat{\Sigma}(u_n) C_{n,k_1,k_2} . \]


\section{Some numerical results}                          
This section presents some results of a simulation study and a real data example. 
For a sample size $n$, the statistic $\widehat{Q}_{n}$ is computed  with $u_n=[\left(\log(n)\right)^{5/2}]$ and $v_n=[\left(\log(n)\right)^{2}]$ (see also Remark 1 in \nocite{Kengne2012}  Kengne (2012)).
The empirical levels and powers are obtained after 200 replications at the nominal level $\alpha=0.05$.

\subsection{Simulation study} 
We consider the following models:

\medskip

(i) ARMA(1,1) processes:
       \begin{equation}\label{arma_model}
       X_t = \alpha_0^* + \alpha_1^* X_{t-1} + \xi_t + \beta_1^* \xi_{t-1} ~ \text{ for all } t \in \Z.
       \end{equation} 
      The parameter of the model is 
      $\theta^*=(\alpha_0^*, \alpha_1^*, \beta_1^*) \in \Theta$, where $\Theta$ is a compact subset of $\R^3$ such as: for all $\theta=(\alpha_0, \alpha_1, \beta_1) \in \Theta$, $|\alpha_1| +  |\beta_1| < 1$.
      Since we can write for all $t \in \Z$,
      \[ X_t= \dfrac{\alpha_0^*}{ 1 + \beta_1^*}   + (\alpha_1^* + \beta_1^*)\big( X_{t-1} +  \sum_{k\geq 2}  (-\beta_1^*)^{k-1} X_{t-k} \big) + \xi_{t},\] 
the  model (\ref{arma_model}) belongs to the class $\mathcal{AC}_{\Z}(M_{\theta^*},f_{\theta^*})$ with $f_\theta (x_1,\cdots) = \dfrac{\alpha_0}{ 1 + \beta_1}   + (\alpha_1 + \beta_1)\big( x_1 +  \sum_{k\geq 2}  (-\beta_1)^{k-1} x_k \big)$ and $M_\theta \equiv 1$ for all $\theta =(\alpha_0, \alpha_1, \beta_1) \in \Theta$.
For this model, the Lipschitz-type conditions A$_i(\Psi_\theta, \Theta)$ ($i=0,1,2$) as well as D$(\Theta)$ are automatically satisfied. 
 Moreover, if $\xi_0$ is a non-degenerate random variable, then the assumptions Id($\Theta$) and  Var($\Theta$)  hold; and for any $r\geq 1$ such
that  $\E \xi_0 ^r < \infty$, $\Theta(r)= \Theta$.
In the sequel, we deal with an ARMA(1,1) with a non zero mean ($\theta^*=(\alpha_0^*, \alpha_1^*, \beta_1^*)$), an ARMA(1,1) with mean zero ($\theta^*=(\alpha_1^*, \beta_1^*)$) and an AR(1) with a non zero mean ($\theta^*=(\alpha_0^*, \alpha_1^*)$). 

\medskip

\begin{figure}[h!]
\begin{center}
\includegraphics[height=10.05cm, width=17.5cm]{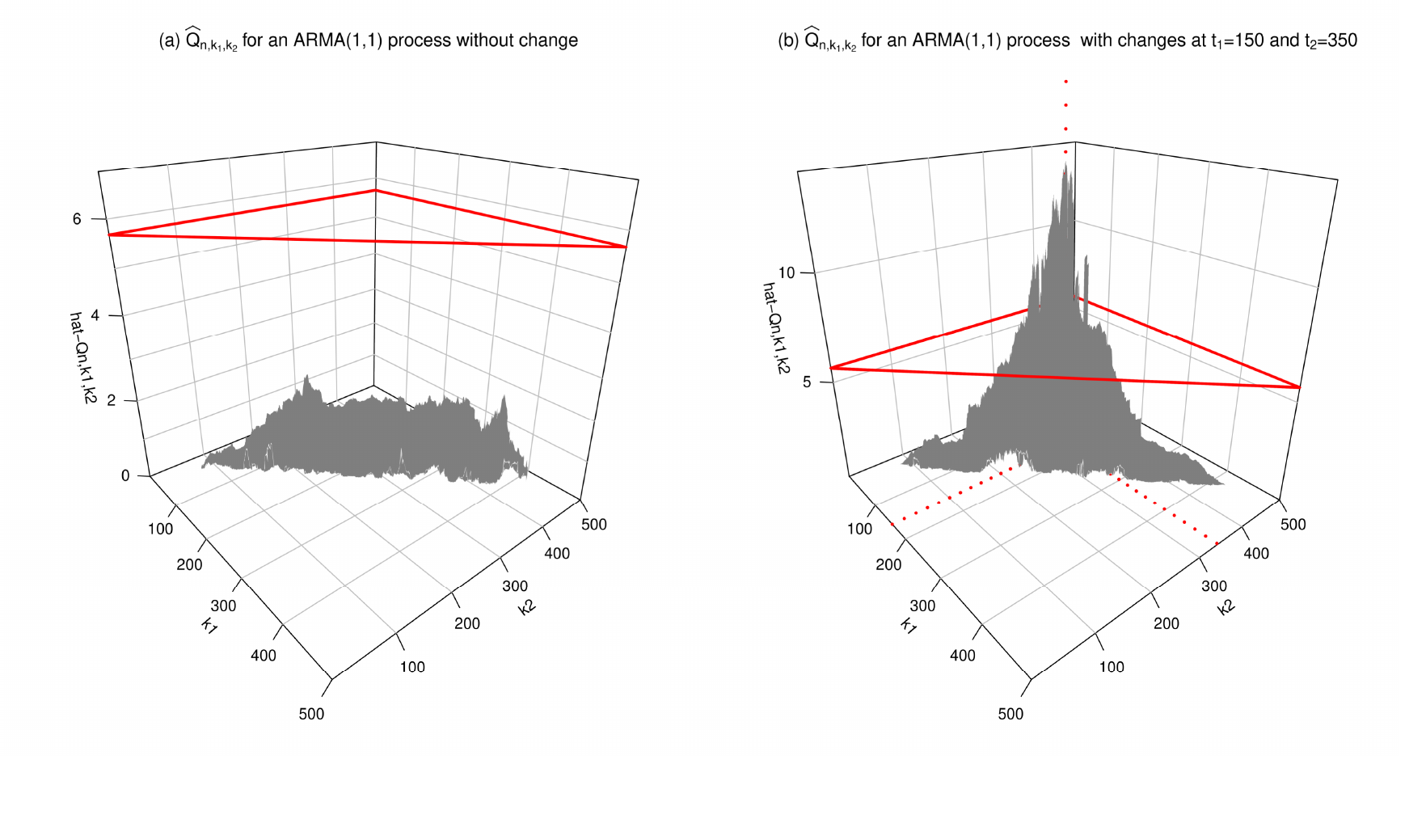} 
\end{center}
\vspace{-1.5cm}
\caption{\small \it 
 The statistics $\widehat{Q}_{n,k_1,k_2}$ for trajectories of an ARMA(1,1) with mean zero: (a) a scenario without change, where the true parameter $\theta^*_0=(-0.4,-0.25)$ is constant; (b)  a scenario under the epidemic alternative where the parameter $\theta^*_1=(-0.4,-0.25)$ changes to $\theta^*_2=(-0.4,0.1)$ at $t_1=150$ and  reverts back to $\theta^*_1$ at $t_2=350$. }
\label{Graphe_sim_ARMA}
\end{figure}

We consider the change-point test with an epidemic alternative where the parameter of the model is $\theta^*_0$ under H$_0$, and $\theta^*_1$, $\theta^*_2$  under H$_1$. 
Firstly, two trajectories of an ARMA(1,1) with mean zero are generated:  a trajectory under H$_0$ with $\theta^*_0=(-0.4,-0.25)$ and a trajectory under H$_1$ with breaks at $(t^*_1, t^*_2) = (150,350)$, $\theta^*_1=(-0.4,-0.25)$, $\theta^*_2=(-0.4,0.1)$. 
 Figure \ref{Graphe_sim_ARMA} displays the statistic $\widehat{Q}_{n,k_1,k_2}$.  
One can see that, for the scenario without change, the values of this statistic are below the horizontal triangle which represents the limit of the critical region (see Figure \ref{Graphe_sim_ARMA}(a)). 
 Under the epidemic alternative,
 $\max_{(k_1, k_2) \in \mathcal T_n} \widehat{Q}_{n,k_1,k_2}$ is greater than the critical value of the test and is reached around the points where the changes occur (see the dotted lines in Figure \ref{Graphe_sim_ARMA}(b)).

 \medskip
 
(ii) GARCH(1,1) processes:  
 \begin{equation}\label{garch_model}
     X_t= \sigma_t \xi_t ~~ \text{ with} ~~ \sigma^2_t= \alpha^*_0 + \alpha^*_1 X^2_{t-1} + \beta^*_1\sigma^2_{t-1},
\end{equation}
the parameter $ \theta^* =(\alpha^*_0,\alpha^*_1,\beta^*_1) \in \Theta$, a compact subset of $]0,\infty[\times[0,\infty[^2$ such as:
for all $\theta=(\alpha_0, \alpha_1, \beta_1) \in \Theta$, $\alpha_1 +\beta_1 < 1$. For all $t \in \Z$, we get
\[ X_t = \xi_t \sqrt{ \alpha^*_0/(1-\beta^*_1) + \alpha^*_1 \sum\limits_{ k\geq 1}^{} (\beta^*_1)^{k-1} X^2_{t-k} }  .\]
Therefore, the  model (\ref{garch_model}) belongs to the class $\mathcal{AC}_{\Z}(M_{\theta^*},f_{\theta^*})$ with $M_\theta (x_1,\cdots) = \sqrt{ \frac{\alpha_0}{ 1 - \beta_1}   + \alpha_1 \sum_{k\geq 1}  \beta_1^{k-1} x_k^2 }$ and $f_\theta \equiv 0$  for all $\theta =(\alpha_0, \alpha_1, \beta_1) \in \Theta$. 
The Lipschitz-type conditions A$_i(\Psi_\theta, \Theta)$ ($i=0,1,2$) hold automatically and D$(\Theta)$ is satisfied with 
 $\underline{h} = \underset{ \theta=(\alpha_0, \alpha_1, \beta_1) \in \Theta  } \inf  \sqrt{ \frac{\alpha_0}{ 1 - \beta_1} } > 0 $. 
 In addition, if $\xi_0$ is a non-degenerate random variable, then the assumptions Id($\Theta$) and  Var($\Theta$)  hold; and for any $r\geq 1$ such
that  $\E \xi_0 ^r < \infty$, $\Theta(r)= \{ \theta=(\alpha_0, \alpha_1, \beta_1) \in \Theta ; ~ \| \xi_0\|_r^2 (\alpha_1 + \beta_1) <1 \}$.
In the sequel, we consider a GARCH(1,1)  ($\theta^*=(\alpha_0^*, \alpha_1^*, \beta_1^*)$) and an ARCH(1) ($\theta^*=(\alpha_0^*, \alpha_1^*)$). 

\medskip

For both the ARMA and GARCH model, we carry out the change-point test with an epidemic alternative where the parameter of the model is $\theta^*_0$ under H$_0$, and $\theta^*_1$, $\theta^*_2$  under H$_1$ with change-points at $(t^*_1, t^*_2) = (0.3n, 0.7n)$ for sample size $n=500, 1000$. The empirical levels and powers are displayed in Table \ref{tab}. The AR$(1)$ example is related to the real data application, see subsection \ref{sub_sec_Real_data}.
The results in this table show that, the empirical level approaching the nominal one when $n$ increases and the empirical power increases with $n$ and is overall close to one when $n = 1000$. These findings are consistent with
the asymptotic results of Theorems \ref{th1} and \ref{th2}.
\begin{table}[h!]
\centering
\scriptsize
\caption{\it Empirical sizes and powers  for the epidemic change-point detection in the models (\ref{arma_model}) and (\ref{garch_model}).}
\label{tab}
\vspace{.2cm}
\hspace*{-.3cm}
\begin{tabular}{llllcc}
\hline
\rule[0cm]{0cm}{.5cm}
&&&&$n=500$&$n=1000$\\
\hline
\rule[0cm]{0cm}{.4cm}
\multirow{3}*{\footnotesize AR(1)}
 &Empirical levels:  &$\theta_0^*=(813,0.3) $       &&0.025&0.035\\
                 
& && &&\\
 
\rule[0cm]{0cm}{.1cm}
&Empirical powers:  &$\theta_1^*=(813,0.3);$&$\theta_2^*=(933,0.24);$ &1.000&1.000\\
                  
          & &&&&\\
   \hdashline[3pt/3pt]  
   & &&&&\\               
\multirow{3}*{\footnotesize ARMA(1,1) with zero mean}
 &Empirical levels:  &$\theta_0^*=(-0.4,-0.25) $       &&0.045&0.050\\
                 
& && &&\\
 
\rule[0cm]{0cm}{.1cm}
&Empirical powers:  &$\theta_1^*=(-0.4,-0.25);$&$\theta_2^*=(-0.4,0.1);$ &0.760&0.990\\

          & &&&&\\
   \hdashline[3pt/3pt]  
   & &&&&\\               
\multirow{3}*{\footnotesize ARMA(1,1) with non zero mean}
 &Empirical levels:  &$\theta_0^*=(1,0.15,0.2) $       &&0.065&0.060\\
                 
& && &&\\
 
\rule[0cm]{0cm}{.1cm}
&Empirical powers:  &$\theta_1^*=(1,0.15,0.2);$&$\theta_2^*=(1,0.5,0.2);$ &0.925&1.000\\

          & &&&&\\
   \hdashline[3pt/3pt]  
   & &&&&\\               
\multirow{3}*{\footnotesize ARCH(1)}
 &Empirical levels:  &$\theta_0^*=(0.6,0.4) $       &&0.035&0.045\\
                 
& && &&\\
 
\rule[0cm]{0cm}{.1cm}
&Empirical powers:  &$\theta_1^*=(0.6,0.4);$&$\theta_2^*=(0.2,0.4);$ &0.910&0.995\\

          & &&&&\\
   \hdashline[3pt/3pt]  
   & &&&&\\               
\multirow{3}*{\footnotesize GARCH(1,1) }
 &Empirical levels:  &$\theta_0^*=(0.15,0.3,0.25) $       &&0.080&0.060\\
                 
& && &&\\
 
\rule[0cm]{0cm}{.1cm}
&Empirical powers:  &$\theta_1^*=(0.15,0.3,0.25);$&$\theta_2^*=(0.15,0.3,0.55);$ &0.730&0.920\\
 
 \Xhline{.75pt}
\end{tabular}
\end{table}

 \subsection{Real data example}\label{sub_sec_Real_data}         
\begin{figure}[h!]
\begin{center}
\includegraphics[height=7.5cm, width=18cm]{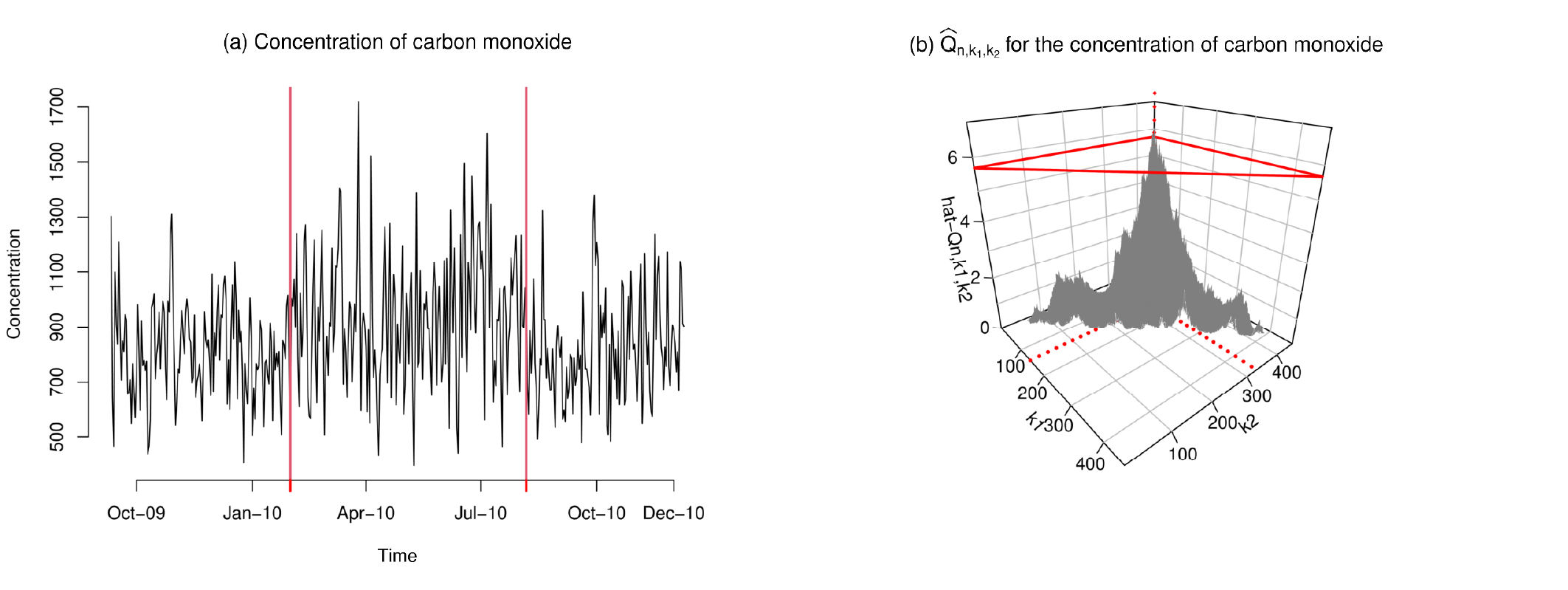} 
\end{center}
\vspace{-.8cm}
\caption{\small \it Plot of $\widehat{Q}_{n,k_1,k_2}$ for the epidemic change-point detection applied to 
the daily averages concentrations for carbon monoxide in the Vit\'oria metropolitan area, between September 11, 2009 and December 09, 2010 with an AR(1) process.
The vertical lines in (a) are the estimated breakpoints. 
The horizontal triangle in (b) represents the limit of the critical region of the test, whereas 
the dotted lines show  the point where the maximum of $\widehat{Q}_{n,k_1,k_2}$ is reached.}
\label{Graphe_Application}
\end{figure}

 We consider the daily concentrations of carbon monoxide in the Vit\'oria metropolitan area.
    These daily levels are obtained from the State Environment and Water Resources Institute, where the data were collected at eight monitoring stations.
     There are $455$ available observations that represent the average concentrations from September 11, 2009 through December 09, 2010 (see Figure \ref{Graphe_Application}(a)). 
The data are a part of a large dataset (available at https://rss.onlinelibrary.wiley.com/pb-assets/hub-assets/rss/Datasets/RSSC\%2067.2/C1239deSouza-1531120585220.zip) which were analyzed by \nocite{Souza2018} Souza {\it et al.} (2018) to quantify the association between respiratory disease and air pollution concentrations.

To test the presence of an epidemic change in this series, we apply our detection procedure with the ARMA($p,q$) model. We have applied the test with several values of $p$ and $q$; and the results after change-point detection show a preference (in the sense of AIC and BIC) for an AR(1). 
%
 Figure \ref{Graphe_Application}(b) displays the values of  $\widehat{Q}_{n,k_1,k_2}$ for all $(k_1,k_2) \in \mathcal T_n$.  
  The critical value on nominal level $\alpha=5\%$ is $c_{d,\alpha}=5.69$ and the resulting test statistic is $\widehat{Q}_{n}=6.61$; which implies that the null hypothesis H$_0$ is rejected (i.e., an epidemic change-point is detected). 
 The vector of the break-points estimated is $\widehat{ \underline t}_n=(143,330)$; i.e, the point where the peak in the graph is reached (see Figure \ref{Graphe_Application}(b)).   
 The locations of the changes correspond to the dates January 31 and August 06, 2010. This corresponds to the period where the winds are weaker and the austral winter; these meteorological factors are noticeable to increase the concentration of the carbon monoxide. 
 The estimated model on each regime is given by:
 \begin{equation}\label{Estim_real_data}
X_{t}=\left\{
\begin{array}{ll}
  \underset{(19.72)}{813.39} + \underset{(0.08)}{0.309}  X_{t-1} +\xi_t~\text{ for }~t\leq 143,  
 \\
\rule[0cm]{0cm}{.6cm}
\underset{(22.43)}{933.27} + \underset{(0.07)}{0.240}  X_{t-1}+\xi_t~ ~\text{ for }~ 144 \leq t  \leq 330   , 
  \\
\rule[0cm]{0cm}{.6cm}
\underset{(22.82)}{822.83} + \underset{(0.09)}{0.293}  X_{t-1}+\xi_t~ ~\text{ for }~  t  \geq 331  , 
\end{array}
\right.
\end{equation}
where in parentheses are the standard errors of the estimators.  
From (\ref{Estim_real_data}), one remark that, the parameter of the first regime is close to that of the third regime; which strengthens the hypothesis of the existence of an epidemic change-point.

 \section{Proofs of the main results} 
  To simplify the expressions, in this section, we will use the conditional Gaussian quasi-log-likelihood up to multiplication by 1/2, given by 
  $L(T,\theta):=-\sum\limits_{t \in T}  q_t(\theta)$ and $\widehat  L(T,\theta):=-\sum\limits_{t \in T}   \widehat q_t(\theta)$.
%

 \subsection{Proof of Theorem \ref{th1}}
 Let $\Sigma :={F}^{-1}\,G\,{F}^{-1}$, where $F$ and $G$ are the matrices defined in (\ref{FG}). 
 Define the statistic
  \[ 
  Q_n=\max_{(k_1, k_2) \in \mathcal T_n} Q_{n,k_1,k_2}
~\text{ with }~
 Q_{n,k_1,k_2}=C_{n,k_1,k_2}' \Sigma \,C_{n,k_1,k_2}.
  \]
  Consider the following lemma; we can go along similar lines as in the proof of Lemma 6.3 in Diop and Kengne (2021) to show the part (i).
   The part (ii) is established in Bardet and Wintenberger (2009). 
 \begin{lem}\label{Lem1} 
 Suppose that the assumptions of Theorem \ref{th1} hold. 
Then, 
\begin{enumerate}
\rm
\item [(i)] $\max_{(k_1,k_2) \in \mathcal T_n}\big|\widehat Q_{n,k_1,k_2}-Q_{n,k_1,k_2}\big|=o_P(1)$;

\rm
  \item [(ii)]  \it $\left(\frac{\partial}{\partial \theta} q_t(\theta^*_0),\mathcal{F}_{t}\right)_{t \in \mathbb{Z}}$ is a stationary ergodic, square integrable martingale difference sequence with covariance matrix $G$.   
\end{enumerate}
 \end{lem}  

\medskip
\noindent
Let two integers $k,k^{\prime} \in[1,n]$, $\bar{\theta} \in \Theta$ and $i \in \{1,2,\cdots,d\}$. Applying the mean value theorem to $\theta \mapsto \frac{\partial}{\partial \theta_i} L({T_{k,k^{\prime}}},\theta)$, there exists $\theta_{n,i}$ between $\bar{\theta}$ and $\theta^*_0$ such that
\[ 
\frac{\partial}{\partial \theta_i} L({T_{k,k^{\prime}}},\bar{\theta})=\frac{\partial}{\partial \theta_i} L({T_{k,k^{\prime}}},\theta^*_0) +\frac{\partial^{2}}{\partial \theta\partial \theta_i} L({T_{k,k^{\prime}}},\theta_{n,i})(\bar{\theta}-\theta^*_0);
\]
i.e,
\begin{eqnarray}\label{Eq_Taylor}
(k^{\prime}-k+1)F_n({T_{k,k^{\prime}}},\bar{\theta}) (\bar{\theta}-\theta^*_0)=\frac{\partial}{\partial \theta} L({T_{k,k^{\prime}}},\theta^*_0)-\frac{\partial}{\partial \theta} L({T_{k,k^{\prime}}},\bar{\theta})
\end{eqnarray}
with
\begin{equation}\label{test_Def_F_n}
F_n({T_{k,k^{\prime}}},\bar{\theta})= - \frac{1}{(k^{\prime}-k+1)}\frac{\partial^{2}}{\partial \theta\partial \theta_i} L({T_{k,k^{\prime}}},\theta_{n,i})_{1\leq i \leq d}.
\end{equation}
 The following lemma will be useful in the sequel.
\begin{lem}\label{Lem2}
Suppose that the assumptions of Theorem \ref{th1} hold.
If $(j_n)_{n\geq1}$  and $(k_n)_{n\geq1}$ are two  integer valued sequences such that $j_n < k_n$, $k_n \limiten \infty$ and $k_n - j_n \limiten \infty $, then
\begin{enumerate}
\rm
\item [(i)] 
 $\frac{1}{\sqrt{k_n-j_n}}\big\| \frac{\partial}{\partial \theta}\widehat{L}(T_{j_n,k_n},\theta)  -\frac{\partial}{\partial \theta}L(T_{j_n,k_n},\theta)\big\|_{\Theta}  \limiteproban 0$;
 
\rm
\item [(ii)]
$
F_n(T_{j_n,k_n},\widehat{\theta}(T_{j_n,k_n})) \limitepsn F$.
\end{enumerate}
\end{lem}

\medskip
\noindent 
{\bf Proof.}
\begin{enumerate}
\item [(i)] See the proof of Theorem 2 in Bardet and Wintenberger (2009).

\item [(ii)]
Applying $(\ref{test_Def_F_n})$ with $\bar{\theta} = \widehat{\theta}(T_{j_n,k_n})$, we obtain
\begin{align*}
F_n(T_{j_n,k_n},\widehat{\theta}(T_{j_n,k_n}))
&=  \Big(\frac{1}{k_n-j_n+1}\frac{\partial^{2}}{\partial \theta\partial \theta_i} L(T_{j_n,k_n},\theta_{n,i})\Big)_{1\leq i \leq d}
=
\frac{1}{k_n-j_n+1}\Big(\sum_{t \in T_{j_n,k_n}}\frac{\partial^{2} q_t(\theta_{n,i})}{\partial \theta\partial \theta_i} \Big)_{1\leq i \leq d},
\end{align*}
where $\theta_{n,i}$ belongs between $\widehat{\theta}(T_{j_n,k_n})$ and $\theta^*_0$.  
Since $\widehat{\theta}(T_{j_n,k_n})  \limitepsn \theta^*_0$, $~ \theta_{n,i} ~ \limitepsn \theta^*_0$ for any  $i=1,\cdots,d$ (from the consistency of the QMLE)
 and that $F=\E\big[\frac{\partial^{2} q_0(\theta^*_0)}{\partial \theta\partial\theta'} \big]$ exists (see Bardet and Wintenberger (2009)), by the uniform strong law of large numbers, for any $i =1,\cdots,d$, we get 
\begin{align*}
 & \Big \| \frac{1}{k_n-j_n+1}\sum_{t \in T_{j_n,k_n}}\frac{\partial^{2}  q_t(\theta_{n,i})}{\partial \theta\partial \theta_i} 
-  \E\Big[\frac{\partial^{2} q_0(\theta^*_0)}{\partial \theta\partial\theta_i}\Big] \Big \| \\
&\leq \Big \| \frac{1}{k_n-j_n+1}\sum_{t \in T_{j_n,k_n}}\frac{\partial^{2}  q_t(\theta_{n,i})}{\partial \theta\partial \theta_i} 
-  \E\Big[\frac{\partial^{2} q_0(\theta_{n,i})}{\partial \theta\partial\theta_i}\Big] \Big \| 
+  \Big \| \E\Big[\frac{\partial^{2} q_0(\theta_{n,i})}{\partial \theta\partial\theta_i}\Big] 
-  \E\Big[\frac{\partial^{2} q_0(\theta^*_0)}{\partial \theta\partial\theta_i}\Big] \Big \|  \\
&\leq \Big \| \frac{1}{k_n-j_n+1}\sum_{t \in T_{j_n,k_n}}\frac{\partial^{2}  q_t(\theta)}{\partial \theta\partial \theta_i} 
-  \E\Big[\frac{\partial^{2} q_0(\theta)}{\partial \theta\partial\theta_i}\Big] \Big \|_\Theta 
+  o(1) = o(1) + o(1)= o(1).
\end{align*}
This establishes the lemma.
\end{enumerate}
\begin{flushright}
$\blacksquare$ 
\end{flushright}

\noindent  
Now, let us show that
\begin{equation}\label{Cond_proof_th1}
Q_n \limiteloin \sup_{0\leq \tau_1<\tau_2 \leq 1}\left\|W_d(\tau_1)-W_d(\tau_2)\right\|^{2}.
\end{equation}
\medskip
Remark that, $Q_{n,k_1,k_2}$ can be rewritten as 
\[
Q_{n,k_1,k_2}=\big\|G^{-1/2}F \cdot C_{n,k_1,k_2}\big\|^2
\]
with
\[
C_{n,k_1,k_2}= \frac{k_2-k_1}{n^{3/2}} \left[ (n-k_2)\left(\widehat{\theta}(T_{k_1+1,k_2})-\widehat{\theta}(T_{k_2+1,n})\right) -  k_1\left( \widehat{\theta}(T_{1,k_1}) -\widehat{\theta}(T_{k_1+1,k_2})\right) \right].
\]
 Let $(k_1,k_2) \in \mathcal T_n$. Applying (\ref{Eq_Taylor}) with $\bar{\theta}=\widehat{\theta}(T_{k_1+1,k_2})$ and $T_{k,k^\prime}=T_{k_1+1,k_2}$, we obtain
\begin{eqnarray}\label{Eq_T_(k_1+1,k_2)}
F_n(T_{k_1+1,k_2},\widehat{\theta}(T_{k_1+1,k_2})) \cdot (\widehat{\theta}(T_{k_1+1,k_2})-\theta_0^*)=\frac{1}{k_2-k_1}\Big(\frac{\partial}{\partial \theta} L(T_{k_1+1,k_2},\theta_0^*)-\frac{\partial}{\partial \theta} L(T_{k_1+1,k_2},\widehat{\theta}(T_{k_1+1,k_2}))\Big).
\end{eqnarray}
With $\bar{\theta}=\widehat{\theta}(T_{k_2+1,n})$ and $T_{k,k^\prime}=T_{k_2+1,n}$, (\ref{Eq_Taylor}) implies  
\begin{eqnarray}\label{Eq_T_(k_2+1,n)}
F_n(T_{k_2+1,n},\widehat{\theta}(T_{k_2+1,n})) \cdot (\widehat{\theta}(T_{k_2+1,n})-\theta_0^*)=\frac{1}{n-k_2}
\Big(\frac{\partial}{\partial \theta} L(T_{k_2+1,n},\theta_0^*)-\frac{\partial}{\partial \theta} L(T_{k_2+1,n},\widehat{\theta}(T_{k_2+1,n}))\Big).
\end{eqnarray}
Moreover, as $n\rightarrow +\infty$, from the asymptotic normality  of the QMLE (see Bardet and Wintenberger (2009)) and Lemma \ref{Lem2}(ii), we have 
\begin{equation}\label{cov_OP_H0}
   \left\{
\begin{array}{l} 
 \big\|\sqrt{k_1} \left(\widehat{\theta}(T_{1,k_1}) -\theta^*_0\right)\big\|=O_P(1), ~~
 \big\|\sqrt{k_2-k_1} \left(\widehat{\theta}(T_{k_1+1,k_2}) -\theta^*_0\right)\big\|=O_P(1);\\

\rule[0cm]{0cm}{.6cm}
\big\|\sqrt{n-k_2} \left(\widehat{\theta}(T_{k_2+1,n}) -\theta^*_0\right)\big\|=O_P(1);\\

\rule[0cm]{0cm}{.6cm}
\big\|F_n(T_{k_1+1,k_2},\widehat{\theta}(T_{k_1+1,k_2}))-F\big\|=o(1)
~\text{ and }~
\big\|F_n(T_{k_2+1,n},\widehat{\theta}(T_{k_2+1,n}))-F \big\|=o(1).
\end{array} 
\right.
\end{equation}
Then, for $n$ large enough, it holds from (\ref{Eq_T_(k_1+1,k_2)}) that 
\begin{align*}
\sqrt{k_2-k_1}F \cdot \left(\widehat{\theta}(T_{k_1+1,k_2})-\theta^*_0\right)
&=
\frac{1}{\sqrt{k_2-k_1}}\Big(\frac{\partial}{\partial \theta} L(T_{k_1+1,k_2},\theta_0^*)-\frac{\partial}{\partial \theta} L(T_{k_1+1,k_2},\widehat{\theta}(T_{k_1+1,k_2}))\Big) 
\\
& \hspace{1cm} -\sqrt{k_2-k_1}\big(\left(F_n(T_{k_1+1,k_2},\widehat{\theta}(T_{k_1+1,k_2}))-F\right)\left(\widehat{\theta}(T_{k_1+1,k_2})-\theta_0\right)\big)\\
&=
\frac{1}{\sqrt{k_2-k_1}}\Big(\frac{\partial}{\partial \theta} L(T_{k_1+1,k_2},\theta_0^*)-\frac{\partial}{\partial \theta} L(T_{k_1+1,k_2},\widehat{\theta}(T_{k_1+1,k_2}))\Big) +o_P(1)\\
&=
 \frac{1}{\sqrt{k_2-k_1}}\Big(\frac{\partial}{\partial \theta} L(T_{k_1+1,k_2},\theta_0^*)- \frac{\partial}{\partial \theta} \widehat L(T_{k_1+1,k_2},\widehat{\theta}(T_{k_1+1,k_2}))\Big)+o_P(1)\\
& \hspace{1cm} 
+\frac{1}{\sqrt{k_2-k_1}}\Big(\frac{\partial}{\partial \theta} \widehat L(T_{k_1+1,k_2},\widehat{\theta}(T_{k_1+1,k_2}))-\frac{\partial}{\partial \theta} L(T_{k_1+1,k_2},\widehat{\theta}(T_{k_1+1,k_2}))\Big) \\
&=
\frac{1}{\sqrt{k_2-k_1}}\Big(\frac{\partial}{\partial \theta} L(T_{k_1+1,k_2},\theta_0^*)-\frac{\partial}{\partial \theta} \widehat L(T_{k_1+1,k_2},\widehat{\theta}(T_{k_1+1,k_2}))\Big) +o_P(1),
\end{align*}
where the last equality is obtained from Lemma \ref{Lem2}(i). 
This is equivalent to
\begin{equation}\label{Eq_a}
F\cdot \left(\widehat{\theta}(T_{k_1+1,k_2})-\theta_0^* \right)=\frac{1}{k_2-k_1}\Big(\frac{\partial}{\partial \theta} L(T_{k_1+1,k_2},\theta_0^*)-\frac{\partial}{\partial \theta} \widehat L(T_{k_1+1,k_2},\widehat{\theta}(T_{k_1+1,k_2}))\Big) +o_P\Big(\frac{1}{\sqrt{k_2-k_1}}\Big).
\end{equation}
 For $n$ large enough, $\widehat{\theta}(T_{k_1+1,k_2})$ is an interior point of $\Theta$ and we have $ \frac{\partial}{\partial \theta} \widehat L(T_{k_1+1,k_2},\widehat{\theta}(T_{k_1+1,k_2}))=0$.
Thus, from (\ref{Eq_a}), we obtain
\begin{equation}\label{Eq_a_bis}
F\cdot \left(\widehat{\theta}(T_{k_1+1,k_2})-\theta_0^*\right)
 =\frac{1}{k_2-k_1}\frac{\partial}{\partial \theta} L(T_{k_1+1,k_2},\theta_0^*)+o_P\Big(\frac{1}{\sqrt{k_2-k_1}}\Big).
\end{equation}
Similarly, using (\ref{Eq_T_(k_2+1,n)}), we also obtain
\begin{equation}\label{Eq_a_bisbis}
F\cdot \left(\widehat{\theta}(T_{k_2+1,n})-\theta_0^*\right)=\frac{1}{n-k_2}\frac{\partial}{\partial \theta} L(T_{k_2+1,n},\theta_0^*)+o_P\Big(\frac{1}{\sqrt{n-k_2}}\Big).
\end{equation}
The subtraction of  (\ref{Eq_a_bis}) and (\ref{Eq_a_bisbis})  gives
\begin{multline*}
F\cdot \left(\widehat{\theta}(T_{k_1+1,k_2})-\widehat{\theta}(T_{k_2+1,n})\right)
=
\frac{1}{k_2-k_1}\frac{\partial}{\partial \theta} L(T_{k_1+1,k_2},\theta_0^*)-\frac{1}{n-k_2}\frac{\partial}{\partial \theta} L(T_{k_2+1,n},\theta_0^*)\\
+o_P\Big(\frac{1}{\sqrt{k_2-k_1}}+\frac{1}{\sqrt{n-k_2}}\Big);
\end{multline*}
i.e.,
\begin{multline}\label{part1_C_n,k1,k2}
\frac{(k_2-k_1)(n-k_2)}{n^{3/2}} F \cdot \left(\widehat{\theta}(T_{k_1+1,k_2})-\widehat{\theta}(T_{k_2+1,n})\right)=\\
\frac{1}{n^{3/2}}
\Big[
(n-k_2)\frac{\partial}{\partial \theta} L(T_{k_1+1,k_2},\theta_0^*)-(k_2-k_1)\frac{\partial}{\partial \theta} L(T_{k_2+1,n},\theta_0^*)
\Big]+o_P(1).
\end{multline}
By going along similar lines, we can also show that   
\begin{multline}\label{part2_C_n,k1,k2}
\frac{k_1(k_2-k_1)}{n^{3/2}} F \cdot \left(\widehat{\theta}(T_{1,k_1})-\widehat{\theta}(T_{k_1+1,k_2})\right)=\\
\frac{1}{n^{3/2}}
\Big[
(k_2-k_1)\frac{\partial}{\partial \theta} L(T_{1,k_1},\theta_0^*)-k_1\frac{\partial}{\partial \theta} L(T_{k_1+1,k_2},\theta_0^*)
\Big]+o_P(1).
\end{multline}
Combining (\ref{part1_C_n,k1,k2}) and (\ref{part2_C_n,k1,k2}), we get
\begin{align*}
F \cdot C_{n,k_1,k_2}
&=
\frac{1}{n^{3/2}}
\Big[(n-(k_2-k_1))\frac{\partial}{\partial \theta} L(T_{k_1+1,k_2},\theta_0^*)
-(k_2-k_1) \big(\frac{\partial}{\partial \theta} L(T_{k_2+1,n},\theta_0^*)+ \frac{\partial}{\partial \theta} L(T_{1,k_1},\theta_0^*)\big) 
\Big] +o_P(1).  \\
&=
\frac{1}{\sqrt{n}}
\Big[
\frac{\partial}{\partial \theta} L(T_{k_1+1,k_2},\theta_0^*)
-\frac{(k_2-k_1)}{n}L(T_{1,n}) 
\Big] +o_P(1)\\
&=
\frac{1}{\sqrt{n}}
\Big[
\frac{\partial}{\partial \theta} L(T_{1,k_2},\theta_0^*)- \frac{\partial}{\partial \theta} L(T_{1,k_1},\theta_0^*)
-\frac{(k_2-k_1)}{n}L(T_{1,n}) 
\Big] +o_P(1) \\
&=
\frac{1}{\sqrt{n}}
\Big[
\big(\frac{\partial}{\partial \theta} L(T_{1,k_2},\theta_0^*) -\frac{k_2}{n}L(T_{1,n})\big)- 
\big(\frac{\partial}{\partial \theta} L(T_{1,k_1},\theta_0^*)
-\frac{k_1}{n}L(T_{1,n})\big) 
\Big] +o_P(1);
\end{align*}
i.e.,
\begin{equation}\label{Eq_b}
G^{-1/2}F  \cdot C_{n,k_1,k_2}
=
\frac{G^{-1/2}}{\sqrt{n}}
\Big[
\big(\frac{\partial}{\partial \theta} L(T_{1,k_2},\theta_0^*) -\frac{k_2}{n}L(T_{1,n})\big)- 
\big(\frac{\partial}{\partial \theta} L(T_{1,k_1},\theta_0^*)
-\frac{k_1}{n}L(T_{1,n})\big) 
\Big] +o_P(1).
\end{equation}
%
%
According to  Lemma \ref{Lem1}(ii), applying the central limit theorem for the sequence $\left(\frac{\partial}{\partial \theta} q_t(\theta_0^*),\mathcal{F}_{t}\right)_{t \in \mathbb{Z}}$ (see \nocite{Billingsley1968} Billingsley (1968)), we obtain
\begin{align*}
\frac{1}{\sqrt{n}}\Big(\frac{\partial}{\partial \theta}L(T_{1,[n \tau_1]},\theta_0^*)-\frac{[n \tau_1]}{n}\frac{\partial}{\partial \theta}L(T_{1,n },\theta_0^*)\Big)
&=
\frac{1}{\sqrt{n}}\Big(\sum_{t=1}^{[n \tau_1]}\frac{\partial}{\partial \theta}q_t(\theta_0^*)-\frac{[n \tau_1]}{n}\sum_{t=1}^{n }\frac{\partial}{\partial \theta}q_t(\theta_0^*)\Big)\\
&~~
\limiteloin B_{G}(\tau_1)-\tau_1 B_{G}(1),
\end{align*}
where $B_{G}$ is a Gaussian process with covariance matrix $\min(s,t)G$. 
Hence, 
\begin{align*}
\frac{G^{-1/2}}{\sqrt{n}}\Big(\frac{\partial}{\partial \theta}L(T_{1,[n \tau_1]},\theta_0^*)-\frac{[n \tau_1]}{n}\frac{\partial}{\partial \theta}L(T_{1,n },\theta_0^*)\Big)
\limiteloin 
B_{d}(\tau_1)-\tau_1 B_{d}(1)=W_d(\tau_1) 
\end{align*}
 in  $D([0,1])$, where $B_d$ is a $d$-dimensional standard motion, and $W_d$ is a $d$-dimensional Brownian bridge.\\
 Similarly, we get 
  \begin{align*}
\frac{G^{-1/2}}{\sqrt{n}}\Big(\frac{\partial}{\partial \theta}L(T_{1,[n \tau_2]},\theta_0^*)-\frac{[n \tau_2]}{n}\frac{\partial}{\partial \theta}L(T_{1,n },\theta_0^*)\Big)
\limiteloin 
B_{d}(\tau_2)-\tau_2 B_{d}(1)=W_d(\tau_2) .
\end{align*}
Thus, as $n\rightarrow \infty$, it comes from (\ref{Eq_b}) that 
\begin{align*}
Q_{n,[n\tau_1],[n\tau_2]}&=
\big\| G^{-1/2} F \cdot C_{n,[n\tau_1],[n\tau_2]}\big\|^2
\limiteloin 
\left\|W_d(\tau_1)-W_d(\tau_2)\right\|^2~\text{ in }~ D([0,1]).
\end{align*}
Hence, for $n$ large enough, we have
\[
Q_n=\max_{\underset{k_1 <k_2-v_n}{v_n \leq k_1< k_2 \leq n-v_n}} Q_{n,k_1,k_2}=
\sup_{\frac{v_n}{n}\leq \tau_1 < \tau_2 \leq 1-\frac{v_n}{n}} Q_{n,[n\tau_1],[n\tau_2]} \limiteloin 
\sup_{0\leq \tau_1 < \tau_2 \leq 1} \left\|W_d(\tau_1)-W_d(\tau_2)\right\|^2;\]
which establishes (\ref{Cond_proof_th1}). Using Lemma \ref{Lem1}(i), we can conclude the proof of the theorem.
\begin{flushright}
$\blacksquare$ 
\end{flushright}

 \subsection{Proof of Theorem \ref{th2}}
 
 Under the epidemic alternative, 
 $(Y_1,\cdots,Y_n)$ is a trajectory of $Y=\{Y_t , t\in \Z\}$ which belongs to \\
 $ \mathcal{AC}_{ \{\cdots,-1,0,1,\cdots,t_1^* \} }(M_{\theta_1^*},f_{\theta_1^*}) \bigcap \mathcal{AC}_{ \{t^*_1+1,\cdots t^*_2 \} }(M_{\theta_2^*},f_{\theta_2^*}) \bigcap \mathcal{AC}_{ \{t^*_2+1,\cdots  \} }(M_{\theta_1^*},f_{\theta_1^*}) $,
 $(t^*_1,t^*_2)=([\tau^*_1 n],[\tau^*_2 n])$ (with $0<\tau^*_1<\tau^*_2<1$) and $\theta^{*}_1 \neq \theta^{*}_2$.
%
%

\medskip 
\noindent 
For any $n \in \N$, we have  
\begin{flalign*}
\widehat{Q}_n &= \underset{(k_1,k_2) \in \mathcal T_n}{\max}\widehat{Q}_{n,k_1,k_2} \geq \widehat{Q}_{n,t^*_1,t^*_2}=C_{n,t^*_1,t^*_2}' \widehat{\Sigma}(u_n) C_{n,t^*_1,t^*_2}\\
\text{with}\hspace{1.5cm}&\\
 C_{n,t^*_1,t^*_2}&= \frac{t^*_2-t^*_1}{n^{3/2}} \left[ \big(n-(t^*_2-t^*_1)\big)\widehat{\theta}(T_{t^*_1+1,t^*_2}) - \left( t^*_1 \widehat{\theta}(T_{1,t^*_1}) + (n-t^*_2)\widehat{\theta}(T_{t^*_2+1,n})\right) \right]
&&
\end{flalign*}
and
\begin{multline*}\label{Sigma_un}
\widehat{\Sigma}(u_n)=\frac{1}{3} 
\big[
\widehat F(T_{1,u_n})  \widehat G(T_{1,u_n})^{-1}   \widehat F(T_{1,u_n}) 
+
\widehat F(T_{u_n+1,n-u_n})  \widehat G(T_{u_n+1,n-u_n})^{-1}  \widehat F(T_{u_n+1,n-u_n}) \\
+
\widehat F(T_{n-u_n+1,n})  \widehat G(T_{n-u_n+1,n})^{-1}  \widehat F(T_{n-u_n+1,n})
\big].
\end{multline*}
%
%
We can also write, 
%
\begin{equation*}
C_{n,t^*_1,t^*_2}= - \frac{(t^*_2-t^*_1)(n-(t^*_2-t^*_1))}{n^{3/2}} \left( \widehat{\theta}(T_{1,t^*_1}) -\widehat{\theta}(T_{t^*_1+1,t^*_2}) + \dfrac{n-t^*_2}{n - (t^*_2-t^*_1)} \big(  \widehat{\theta}(T_{t^*_2+1,n}) - \widehat{\theta}(T_{1,t^*_1}) \big) \right).
\end{equation*}
Moreover, by definition, the three matrices in the formula of $\widehat{\Sigma}_n(u_n)$ are positive semi-definite. 
Then, according to the assumption \textbf{B} and for $n$ large enough, we can find a constant $C>0$ such that, it holds $a.s.$
\begin{align}
\widehat{Q}_n &\geq \widehat{Q}_{n,t^*_1,t^*_2} \nonumber\\
& \geq 
\frac{(t^*_2-t^*_1)^2(n-(t^*_2-t^*_1))^2}{n^{3}} \left( \widehat{\theta}(T_{1,t^*_1}) -\widehat{\theta}(T_{t^*_1+1,t^*_2}) +  \dfrac{n-t^*_2}{n - (t^*_2-t^*_1)} \big(  \widehat{\theta}(T_{t^*_2+1,n}) - \widehat{\theta}(T_{1,t^*_1}) \big) \right)'\nonumber\\
& \hspace{2cm} \times \Big[
\widehat F(T_{1,u_n})  \widehat G(T_{1,u_n})^{-1}   \widehat F(T_{1,u_n}) +
\widehat F(T_{n-u_n+1,n})  \widehat G(T_{n-u_n+1,n})^{-1}  \widehat F(T_{n-u_n+1,n})
\Big]\nonumber\\
& \hspace{2.5cm} \times  \left( \widehat{\theta}(T_{1,t^*_1}) -\widehat{\theta}(T_{t^*_1+1,t^*_2}) +  \dfrac{n-t^*_2}{n - (t^*_2-t^*_1)} \big(  \widehat{\theta}(T_{t^*_2+1,n}) - \widehat{\theta}(T_{1,t^*_1}) \big)\right)\nonumber\\
& \geq C \times n
 \left( \widehat{\theta}(T_{1,t^*_1}) -\widehat{\theta}(T_{t^*_1+1,t^*_2}) +  \dfrac{n-t^*_2}{n - (t^*_2-t^*_1)} \big(  \widehat{\theta}(T_{t^*_2+1,n}) - \widehat{\theta}(T_{1,t^*_1}) \big)\right)'\nonumber\\
& \hspace{2cm} \times \Big[
\widehat F(T_{1,u_n})  \widehat G(T_{1,u_n})^{-1}   \widehat F(T_{1,u_n}) +
\widehat F(T_{n-u_n+1,n})  \widehat G(T_{n-u_n+1,n})^{-1}  \widehat F(T_{n-u_n+1,n})
\Big] \nonumber\\
& \hspace{2.5cm}  \label{minor_Q_n} \times  \left( \widehat{\theta}(T_{1,t^*_1}) -\widehat{\theta}(T_{t^*_1+1,t^*_2}) +  \dfrac{n-t^*_2}{n - (t^*_2-t^*_1)} \big(  \widehat{\theta}(T_{t^*_2+1,n}) - \widehat{\theta}(T_{1,t^*_1}) \big) \right).
\end{align}
%
From the asymptotic properties of the QMLE, and the study of the stationary regime approximation developed by Bardet {\it et al.} (2012) (see Proposition 6.1 and Corollary 6.1 of these authors), we deduce 
\begin{equation*}
   \left\{
\begin{array}{l} 
\bullet~
\widehat{\theta}(T_{1,t^{*}})-\widehat{\theta}(T_{t^*_1+1,t^*_2}) \limitepsn \theta^{*}_{1}-\theta^{*}_{2}\neq 0;\\

\rule[0cm]{0cm}{.5cm}
\bullet~
\widehat{\theta}(T_{1,u_n}) \limitepsn \theta^{*}_1,
 ~~~\widehat{\theta}(T_{n-u_n+1,n}) \limitepsn \theta^{*}_1;\\
 
\rule[0cm]{0cm}{.5cm}
\bullet~
\dfrac{n-t^*_2}{n - (t^*_2-t^*_1)} \big(  \widehat{\theta}(T_{t^*_2+1,n}) - \widehat{\theta}(T_{1,t^*_1}) \big) \limitepsn 0;\\

\rule[0cm]{0cm}{.5cm}
\bullet~
\widehat F(T_{1,u_n})  \widehat G(T_{1,u_n})^{-1}   \widehat F(T_{1,u_n}) + 
\widehat F(T_{n-u_n+1,n})  \widehat G(T_{n-u_n+1,n})^{-1}  \widehat F(T_{n-u_n+1,n}) \limitepsn 2\Sigma^{(1)},
\end{array} 
\right.
\end{equation*}
where $\Sigma^{(1)}$ denotes the covariance matrix of the stationary model of the first and third regimes, defined as in (\ref{TLC}) and (\ref{FG}) and computed at $\theta^*_1$. This matrix, which is defined on the stationary regime is positive definite (see  Bardet and Wintenberger (2009)). This implies $\widehat{Q}_{n,t^*_1,t^*_2} \limitepsn +\infty$. Thus, the theorem is obtained from (\ref{minor_Q_n}). 
\begin{flushright}
$\blacksquare$
\end{flushright}
                 
\bibliographystyle{acm}

 \end{document}